\documentclass[10pt]{amsart}

\usepackage{xcolor}

\newcommand{\bc}{}

\usepackage{mathabx}
\input bookman.sty
\boldmath

\usepackage{xcolor}

\usepackage{helvet}

\usepackage{graphics}

\usepackage{graphicx}

\usepackage{amssymb}
\usepackage{amsxtra}
\usepackage{amsmath}
\usepackage{mathrsfs}

\usepackage[arrow, matrix]{xy}
\xyoption{frame}

\theoremstyle{plain}

\newtheorem{theo}{Theorem}[section]
\newtheorem{lemm}[theo]{Lemma}
\newtheorem{prop}[theo]{Proposition}
\newtheorem{coro}[theo]{Corollary}

\theoremstyle{definition}

\newtheorem{defi}[theo]{Definition}
\newtheorem{rema}[theo]{Remark}

\newfont{\rmm}{cmr10 scaled 1000}

\newfont{\itt}{cmsl10 scaled 1000}

\newfont{\rM}{cmr10 scaled 1700}

\newcounter{lemma}[section]

\newcounter{tempcounter}

\newcommand{\lb}{\label}

\newcommand{\rrf}[1]{(\ref{#1})}

\renewcommand{\b}{\beta}
\newcommand{\g}{\gamma}

\newcommand{\e}{\epsilon}

\renewcommand{\t}{\theta}
\renewcommand{\l}{\lambda}

\renewcommand{\r}{\rho}

\newcommand{\s}{\sigma}

\newcommand{\Si}{\Sigma}

\newcommand{\MM}{{\mathcal M}}
\newcommand{\NN}{{\mathcal N}}

\newcommand{\qq}{{\mathbb{Q}}}
\newcommand{\rr}{{\mathbb{R}}}

\newcommand{\zz}{{\mathbb{Z}}}

\newcommand{\RRR}{{\mathbf{R}}}

\newcommand{\ind}{\text{\rm ind\hspace{0.05cm}}}

\renewcommand{\Im}{\text{\rm Im }}

\newcommand{\Int}{\text{\rm Int }}

\newcommand{\bere}{\begin{rema}}
\newcommand{\bede}{\begin{defi}}

\renewcommand{\beth}{\begin{theo}}
\newcommand{\bele}{\begin{lemm}}
\newcommand{\bepr}{\begin{prop}}
\newcommand{\beeq}{\begin{equation}}
\newcommand{\bega}{\begin{gather}}
\newcommand{\begaa}{\begin{gather*}}
\newcommand{\been}{\begin{enumerate}}

\newcommand{\bedee}{\begin{defii}}
\newcommand{\bethh}{\begin{theoo}}
\newcommand{\belee}{\begin{lemmm}}
\newcommand{\beprr}{\begin{propp}}

\newcommand{\beco}{\begin{coro}}

\newcommand{\beal}{\begin{aligned}}

\newcommand{\enre}{\end{rema}}

\newcommand{\enco}{\end{coro}}
\newcommand{\enpr}{\end{prop}}
\newcommand{\enth}{\end{theo}}
\newcommand{\enle}{\end{lemm}}
\newcommand{\enen}{\end{enumerate}}
\newcommand{\enga}{\end{gather}}
\newcommand{\engaa}{\end{gather*}}
\newcommand{\eneq}{\end{equation}}
\newcommand{\enal}{\end{aligned}}

\newcommand{\bq}{\begin{equation}}
\newcommand{\bqq}{\begin{equation*}}

\renewcommand{\leq}{\leqslant}
\renewcommand{\geq}{\geqslant}

\newcommand{\wi}{\widetilde}

\newcommand{\ove}{\overline}

\newcommand{\wh}{\widehat}

\newcommand{\sm}{\setminus}

\newcommand{\sbs}{\subset}

\newcommand{\tens}[1]{\underset{#1}{\otimes}}

\newcommand{\sut}{~such~that~}

\newcommand{\wrt}{with respect to}

\newcommand{\ma}{manifold}
\newcommand{\nei}{neighbourhood}

\newcommand{
\sma}{submanifold}

\newcommand{\Prf}{{\it Proof.\quad}}

\newcommand{\smo}{C^{\infty}}

\newcommand{\chart}{\Phi_p:U_p\to B^n(0,r_p)}
\newcommand{\atlas}{\{\Phi_p:U_p\to B^n(0,r_p)\}_{p\in S(f)}}

\newcommand{\pr}{\partial}

\newcommand{\qs}{\hfill\square}

\newcommand{\arrh}[3]
{
\xymatrix{
{#1} \ar[r]^<<<<{#2}  &{#3}
}
}

\newcommand{\arrr}[1]
{\arrh {}{#1}{}}

\newcommand{\arrto}
{\xymatrix{{} \ar@{|-{>}}[r]  & {} } }

\newcommand{\arrinto}
{\xymatrix{{} \ar@{^{(}->}[r]  & {} } }

\newcommand{\mnk}{\MM\NN(K)}

\newcommand{\ssmk}{S^4\sm S(K)}

\newcommand{\rmf}{regular Morse function}
\newcommand{\df}{diffeomorphism}

\begin{document}

\title
[Circle-valued Morse theory for frame spun knots and surface-links]
{Circle-valued Morse theory for frame spun knots and surface-links}
\author{Hisaaki Endo  and  Andrei Pajitnov}
\address{Department of Mathematics
Tokyo Institute of Technology
2-12-1 Ookayama, Meguro-ku
Tokyo 152-8551
Japan}
\email{endo@math.titech.ac.jp}
\address{Laboratoire Math\'ematiques Jean Leray 
UMR 6629,
Facult\'e des Sciences,
2, rue de la Houssini\`ere,
44072, Nantes, Cedex}                    
\email{andrei.pajitnov@univ-nantes.fr}

\begin{abstract}
 Let $N^k\sbs S^{k+2}$ be a closed oriented
 submanifold, denote its complement by $C(N)= S^{k+2}\sm N$.
 Denote by $\xi\in H^1(C(N))$ the class dual to $N$. 
 The Morse-Novikov number of $C(N)$ is by definition
 the minimal possible number of critical points
 of a regular Morse map $C(N)\to S^1$ belonging to $\xi$.
 In the first part of this paper
 we study the case when $N$ is the twist frame spun knot associated
 to an $m$-knot $K$. We obtain a formula which relates the 
 Morse-Novikov numbers
 of $N$ and $K$ and generalizes the classical results
 of D. Roseman and E.C. Zeeman about fibrations of spun knots.
 In the second part we apply the obtained results to 
 the computation of Morse-Novikov numbers of 
 surface-links in 4-sphere. 
\end{abstract}
\maketitle
\tableofcontents

\section{Introduction}
\label{s:intro}
\subsection{Overview of the article}
\label{su:overv}

Let $N^k\sbs S^{k+2}$ be a closed oriented
 submanifold, let $C(N)= S^{k+2}\sm N$ be its complement.
 The orientation of $N$ determines a cohomology class 
 $\xi\in H^1(C(N))\approx [C(N), S^1]$.
 We say that $N$ is fibred if there is a Morse map
 $f:C(N)\to S^1$ 
 homotopic to $\xi$ 
 which is regular nearby $N$ (see Def \ref{d:reg})
 and has no critical points. 
  In general a Morse map $C(N)\to S^1$ has some critical points, 
the minimal number of these critical points
will be called {\it the Morse-Novikov number of $N$}
and denoted $\MM\NN(C(N))$.
 
 In the first part of this paper we study this invariant 
 in relation with constructions of spinning.
 The classical Artin's spinning construction \cite{Ar}
associates to each knot $K\sbs S^3$ a 2-knot $S(K)\sbs S^4$.
A twisted version of this construction is due to E.C. Zeeman \cite{Z}.
In \cite{Roseman}
 D. Roseman introduced a {\it frame spinning } construction, and G. Friedman 
 \cite{FriedmanAlex}
gave a twisted version of 
generalized 
Roseman's construction to include twisting.

The input data for twist frame spinning  construction is:

\been\lb{page:tfs}
\item[(TFS1)]
A closed \ma~ $M^k\sbs S^{m+k}$ with trivial (and framed ) normal bundle. 
\item[(TFS2)] An $m$-knot $K^m\sbs S^{m+2}$.
\item[(TFS3)]
A smooth map $\l: M\to S^1$.
\enen

To these data one associates an $n$-knot 
$\s(M,K,\l)$, where $n=k+m$  (see Section \ref{s:frame-twist}).
We prove in Section
\ref{s:frame-twist}
the following formula:
\begin{equation}\lb{f:tw-sp}
\MM\NN(C(\s(M,K,\l)))
\leq
\MM\NN(C(K))\cdot \MM\NN(M, [\l])
\end{equation}
(where $\MM\NN(M, [\l])$
is the minimal number of critical points
of a map $M\to S^1$ homotopic to $\l$.)
If $\l$ is null-homotopic, we have 
$$
\MM\NN(C(\s(M,K)))
\leq
\MM\NN(C(K))\cdot \MM(M)
$$
(where $\MM(M)$ is the Morse number of $M$).
In particular, if $K$ is fibred, then the framed spun knot 
$\s(M,K)$ is fibred (A theorem due to D. Roseman \cite{Roseman}).
If in formula 
\rrf{f:tw-sp}
the map $l:M\to S^1$ has no critical points,
then the knot $\s(M,K,\l)$ is fibred, and we recover the classical
result of E.C. Zeeman \cite{Z}: for any knot its twist-spun knot
is fibered.
In section \ref{s:rotation} we discuss a geometric construction
which is related to spinning,
namely rotation of a knot $K^m\sbs S^{m+2}$
around equatorial sphere $\Sigma$ of $S^{m+2}$.
The resulting submanifold $R(K)$ is diffeomorphic to $S^1\times S^m$ 
and is sometimes called
{\it spun torus } of $K$. We prove that 
$$
\MM\NN(R(K))
\leq 2\MM\NN(K) + 2.
$$

Section \ref{s:surf-links} is about Morse-Novikov theory for surface-links.
In Subsection \ref{su:saddle}
 we introduce a related invariant of surface-links, 
 namely the {\it
saddle number $sd(F)$} (Definition \ref{d:saddle})
and prove the formula
\begin{equation}\lb{f:mn-saddle}
\MM\NN(C(F))
\leq
2sd(F) +\chi(F) -2.
\end{equation}
In Subsection 
\ref{su:spun}
we discuss the case of spun knots.
In subsection \ref{su:yosh_links}
we apply the results of 
Sections \ref{s:frame-twist}
and
\ref{s:rotation}
and the formula 
\rrf{f:mn-saddle}
to determine the Morse-Novikov 
numbers of certain surface-links.
In
\cite{Yoshikawa}
K. Yoshikawa
introduced a numerical invariant $ch(F)$ 
of surface-links $F$ and developed a method
that allowed him to enumerate all the (weakly prime) 
surface-links $F$ with $ch(F)\leq 10$.
In Subsection \ref{su:yosh_links}
we compute the Morse-Novikov numbers 
of the majority of the oriented surface-links 
of the Yoshikawa's table.

\subsection{Basic definitions and lower bounds for Morse-Novikov numbers}
\lb{su:basics}

We start with the definition of a regular Morse map.

\begin{defi}\lb{d:reg}
Let $N^k\sbs S^{k+2}$ be a closed oriented
 submanifold. Denote by $\xi\in H^1(C(N))\approx [C(N), S^1]$
 the cohomology class dual to the orientation of $N$.
 A Morse map $f ~:~C(N)\to S^1$ is said to be {\it regular\/}
if there is  an orientation preserving  $\smo$ trivialisation   
\begin{equation}\lb{f:triv}
\Phi:T(N) \to N \times B^2(0,\e)
\end{equation}
of a tubular neighbourhood $T(N)$ of $N$
such that the 
restriction $f|\:\big( T(N)\sm N \big)$
satisfies 
$f\circ \Phi^{-1} (x,z)=  z/|z|.$

An $f$-gradient $v$ of a regular Morse map $f~:~C(N)\to S^1$
will be called {\it regular} if there is a $\smo$ trivialisation   
\rrf{f:triv} 
such that $\Phi^*(v)$ equals $(0,v_0)$
where $v_0$ is the Riemannian gradient of 
the function $z\mapsto  z/|z|$. 

\end{defi}

If $f$ is a Morse map of a manifold to 
$\RRR$ or to $S^1$,
then we denote by $m_p(f)$
the number of critical points of $f$ of index $p$.
The number of all critical points of $f$
is denoted by $m(f)$.

\bede\lb{d:def-MN}
The minimal  number 
$m(f)$ where $f:C(N)\to S^1$ is a 
regular Morse map is called 
{\it the Morse-Novikov number of $N$}
and denoted by $\MM\NN(C(N))$.
 \end{defi}
 
 To obtain lower bounds for numbers $m_p(f)$
 one uses the {\it Novikov homology}.
 Let $L=\zz[t, t^{-1}]$;  denote by  $\wh L=\zz((t))$ and 
$\wh L_\qq=\qq((t))$
the rings of all 
series in one variable $t$ 
with integer (respectively rational) coefficients  and finite negative part. 
Recall that $\wh L$ is a PID, and $\wh L_\qq$ is a field.
Consider the infinite cyclic covering 
$\ove {C(N)} \to C(N)$; the Novikov homology of $C(N)$ is defined as follows: 
$$
\wh H_*(C(N))= H_*(\ove {C(N)})\tens{L}\wh L.
$$
The rank and torsion number of the $\wh L$-module 
$\wh H_*(C(N))$
will be denoted by $\wh b_k(C(N))$, respectively $\wh q_k(C(N))$.
For any  regular Morse function
$f$ there is a Novikov complex 
$\NN_*(f,v)$ 
 over $\wh L$ generated in degree $k$ by critical points of $f$ of index $k$
 and such that 
 $H_*(\NN_*(f,v))\approx \wh H_*(C(N)).$
 (see \cite{CVMT}).
 Therefore we have  the Novikov inequalities
 $$  
 \sum_k \Big(\wh b_k(C(N))+\wh q_k(C(N))+\wh {q}_{k-1}(C(N))\Big)
 \leq
  \MM\NN(C(N)).
 $$
  These inequalities, which are far from being exact in general, are however 
 very useful in particular in the case of surface-links
 (see Section \ref{s:surf-links}).


\section{Twist frame  spun knots}
\lb{s:frame-twist}

We start with a recollection of the twist frame 
 spinning construction following 
\cite{Roseman}, \cite{Fr}, \cite{FriedmanAlex}.
See the input data (TFS1) -- (TFS3) for this construction on the page
\pageref{page:tfs}.
Let $a\in K$.  Removing a small open disk
$D(a)$ from $S^{m+2}$ we obtain an embedded (knotted)
disk $K_0$ in the disk $D^{m+2}\approx S^{m+2}\sm D(a)$.
We identify 
$D^{m+2}$ with the standard Euclidean disk of radius 1 and center $0$
in $\rr^{m+2}$. 
We have the usual diffeomorphism
$$\chi: S^{m+1}\times [0,1[ \arrr {\approx} D^{m+2}\sm \{0\},
\ \ \ (x,t)\mapsto tx.$$
We can assume that $K_0\cap \pr D^{m+2}$
is the standard sphere $S^{m-1}$ in 
$\pr D^{m+2}= S^{m+1}$.
Moreover, we can assume that the intersection of $K_0$ with a \nei~
of $\pr D^{m+2}$ is also standard, that is, 
$$
K_0\cap \chi\big(S^{m+1} \times [1-\e,1]\big)
=
\chi\big(S^{m-1} \times [1-\e,1]\big).
$$
We have a  framing of $M$ in $S^n$;
combining this with the 
standard framing of $S^n$ in $S^{n+2}$ 
we obtain a diffeomorphism 
$$
\Phi: N(M, S^{n+2})\arrr {\approx} M\times D^m\times D^2$$
where $N(M, S^{n+2})$ is  a regular \nei~ of $M$ in $S^{n+2}$.
We can assume that the restriction of $\Phi$ to 
$N(M, S^{n})$ gives a diffeomorphism
$$
\Phi: N(M, S^{n})\arrr {\approx} M\times D^m\times \{0\},$$
induced by the given framing of $M$.
The Euclidean disc $D^{m+2}$ is a subset of $D^m\times D^2$,
so that $K_0\sbs D^m\times D^2$.

For $\t\in S^1$ denote by $R_\t$ the rotation of 
$D^{2}$
around its center.
The disc $D^{m+2}\sbs D^m\times D^2$
is invariant \wrt~ this rotation 
as well as the intersection 
of $K_0$ with a small \nei~ of $\pr D^{m+2}$.
We have 
$\Phi\big(S^n\cap  N(M, S^{n+2})\big)
=
M\times D^m\times \{0\}.$
Let
$$
Z=\{(x,y,z)~|~ (y,z)\in R_{\l(x)}(K_0)\}.$$
This is an $m$-dimensional submanifold of 
$
M\times D^m\times D^2$. We define $\s(M,K,\l)$ as follows
$$
\s(M,K,\l)
=
\Big(S^{n+2} \sm N(M, S^{n+2})\Big)
\cup \Phi^{-1}(Z).
$$
This is the image of an embedded $n$-sphere, knotted in general.

\vskip0.1in

\centerline{\bf Examples and particular cases.}
\been\item
Let $\dim M=0$, so that $M$ is a finite set;
denote by $p$ its cardinality.
Then the $n$-knot $\s(M,K,\l)$ is equivalent to the 
connected sum of $p$ copies of $K$.
\item 
If $M$ is the equatorial circle of the sphere $S^2$,
which is in turn considered as an equatorial sphere of $S^4$,
and $\l(x)=1$, we obtain the classical Artin's construction.
If $\l:S^1\to S^1$ is a map of degree $d$, we obtain the 
Zeeman's twist-spinning construction \cite{Z}.
\item 
If $\l(x)=1$ for all $x\in M$ we obtain the Roseman's
construction of spinning around
the manifold $M$ \cite{Roseman}.
In this case we will denote $\s(M,K,\l)$ by $\s(M,K)$.
\enen

\beth\label{t:mn_tw_spin}
$$
\MM\NN(\s(M,K,\l))\leq \MM\NN(K)\cdot \MM\NN(M, [\l]).
$$
(where $[\l]\in H^1(M,\zz)\approx [M, S^1]$
is the homotopy class of $\l$).
\enth
\Prf 
We will be using the terminology from the above construction
of $\s(M,K,\l)$.
We have the standard fibration
$$\psi_0: 
S^{n+2}\sm S^n \to S^1
$$
obtained from the canonical framing of $S^n$ in $S^{n+2}$.
Observe that the map $\alpha=\psi_0\circ \Phi^{-1}$
is defined by the following formula
$$
\alpha(x,y,z)=\frac z{|z|}.
$$
Let $f:S^{m+2}\sm K\to S^1$ be a Morse map.
The restriction of $f$ to the subset 
$D^{m+2}\sm K_0$ will be denoted by the same letter $f$.
We can assume that the function $f$ equals $\alpha$ 
in a \nei~ of $\pr D^{m+2}=S^{m+1}$.
In particular in a \nei ~ of $\pr D^{m+2}$
we have
$$
f(R_\t(p))=f(p)+\t,
\ \  {\rm for ~ ~} 
p\in S^{m+1}\sm K_0.
$$
Define a function $g$ on 
$M\times D^{m+2}\sm Z$ by the following formula:

\begin{equation}\lb{f:def_g}
g(x,\xi)=f\big(R_{-\l(x)}(\xi)\big) +\l(x),
\end{equation}
(where $x\in M, \ \xi \in D^{m+2}$).
Define a function $\psi$ on the complement $S^{n+2}\sm \s$
by the 
 the following formula:
 \been\item 
 If $p\notin N(M, S^{n+2})$, then $\psi(p)=\psi_0(p)$.
 \item 
 If $p\in N(M, S^{n+2})$, then 
 $\psi(p)=g(\Phi^{-1}(p)).$
 \enen 
We will now prove that if $\l$ is a Morse map
(this can be achieved by a small perturbation of $\l$), then 
$\psi$ is also a Morse map, and the number $m(\psi)$
of its critical points satisfy 
$$
m(\psi) = m(\l)\cdot m(f).
$$
All the critical points of $\psi$ are in 
$N(M, S^{n+2})$. In this domain the function $\psi$ 
is diffeomorphic to $g$, and the count of critical points 
of $g$ is easily achieved with the help of the next lemma.

\bele\lb{l:morse-twisted}
Let $g_1:N_1\to S^1, \ \ g_2:N_2\to S^1$
be Morse functions on manifolds $N_1, N_2$.
Let $F:N_1\times N_2\to N_2$ be a map, such that 
for each $a\in N_2$ the map $x\mapsto F(a,x)$ is a diffeomorphism
$N_2\to N_2$.
Define a function $g:N_1\times N_2\to S^1$ by the following formula:
$$
g(x_1, x_2)= g_1(x_1)+g_2(F(x_1, x_2)).
$$
Then $g$ is a Morse function,
$Crit(g)= Crit(g_1)\times Crit(g_2)$
and for every $a_1\in Crit(g_1), \ a_2\in Crit(g_2)$
we have 
$\ind (a_1, a_2)=\ind (a_1)+\ind(a_2)$.
\enle
\Prf 
Define a function $g_0$ on $N_1\times N_2$ by the following formula
$$
g_0(x_1, x_2)= g_1(x_1)+g_2(x_2).
$$
The conclusions of our Lemma hold obviously
if we replace $g$ by $g_0$ in the statement 
of the Lemma. 
Observe now that the function $g$ is diffeomorphic to $g_0$
via the diffeomorphism
$$
(x_1, x_2)\mapsto (x_1, F(x_1, x_2)).
$$
The lemma follows. $\qs$

\beco\lb{c:class-spin}
Let $K\sbs S^3$ be a classical knot,
denote by $S(K)$ the spun knot of $K$.
Then
\bq
\MM\NN(S(K))\leq 2\MM\NN(K)
\end{equation}
\enco
\Prf
In this case $M=S^1$ and $[\l]=0$. 
We have $\MM\NN(S^1,0)=2$ and the result follows. $\qs$

The classical theorems concerning 
fibrations of spun knots follow from 
Theorem \ref{t:mn_tw_spin}:

\beco\lb{c:roseman}(D. Roseman \ \cite{Roseman})
\ \ \ 
If $K$ is fibred, then $\MM\NN(\s(M,K))$
is fibred.
\enco
\Prf
Since $\MM\NN(K)=0$, Theorem \ref{t:mn_tw_spin}
implies $\MM\NN(\s(M,K))=0$. $\qs$

\beco\lb{c:zeeman}(E.C. Zeeman \ \cite{Z})
\ \ \ 
The $d$-twist spun knot of any classical knot $K$
is fibred for $d\geq 1$.
\enco
\Prf
Consider a great circle $\Sigma$ in $S^2$.
The $d$-twist spun knot of $K$
is by definition the $n+1$-knot $\s(\Sigma, K, \l)$ 
in $S^{3}$
where $\:\Si\to \Si$ is a map of degree $d$.
The assertion follows, since $\MM\NN(S^1,\l)=0$.
$\qs$

\bere\lb{r:general_frame_twist}
The Zeeman's theorem above generalizes immediately to the following statement:
If $\MM\NN(M,\l)=0$, then 
the knot $\s(M,K,\l)$ is fibred for any knot $K$.
\enre



\section{Rotation}
\lb{s:rotation}

Let $\Sigma$ be an equatorial sphere of 
$S^{n+1}$. We can view the sphere 
$S^{n+1}$ as the union of two discs $D_+\cup D^-$
intersecting by $\Sigma$. Consider 
$S^{n+1}$ as the equatorial sphere of 
$S^{n+2}$. The sphere $S^{n+2}$
can be considered  as the result of 
rotation of the disc $D_+$ around its boundary 
$\Sigma$.  We have the (linear orthogonal)
action of $S^1$ on 
$S^{n+2}$, such that $\Sigma$
is the fixed point set of the action,
and the action is free on the rest of 
the sphere $S^{n+2}$.
Let $K^{n-1}$ be an $(n-1)$-knot in
$S^{n+1}$. 
We can assume that $K^{n-1}\sbs \Int D_+$.
Rotation of $K^{n-1}$ around $\Sigma$ 
gives a \sma~ $R(K)$ of codimension 2 in 
$S^{n+2}$. The manifold $R(K)$ is diffeomorphic to $S^1\times K$.
We call this constrution {\it rotation}.
When $\dim K =1$,
the  manifold $R(K)$ 
is sometimes called the
{\it spun torus} of $K$.

In this section we 
relate the Morse-Novikov numbers of 
$R(K)$ with those of $K$. The main aim of this section
is to prove the following theorem.

\beth\lb{t:mn_rotation}
$$
\MM\NN(R(K))
\leq 
2\MM\NN(K) + 2.
$$
\enth
To prove the theorem 
we associate to each given 
\rmf~
$\phi:S^{n+1}\sm K^{n-1} \to S^1$
a \rmf~ 
$R(\phi):S^{n+2}\sm R(K^{n-1}) \to S^1$
such that
$m(R(\phi))=2m(\phi)+2.$

We begin by an outline of this construction
for the simplest case when $n=1$ and $K$ consists 
of two points in $S^2$ (Subsection \ref{su:two_points}).
In Subsection \ref{su:gen_case}
we give a detailed proof of the assertion
of the theorem in  full generality.

\subsection{ Rotation of $S^0$}
\lb{su:two_points}

Let $K^0=\{a,b\} \sbs S^2$. 
The manifold $S^2\sm \{a,b\}$ is fibered over $S^1$, 
and the structure of the level lines of this fibration
is shown on the figure 1 (left).

\begin{figure}
\begin{center}
\includegraphics[height=43mm]{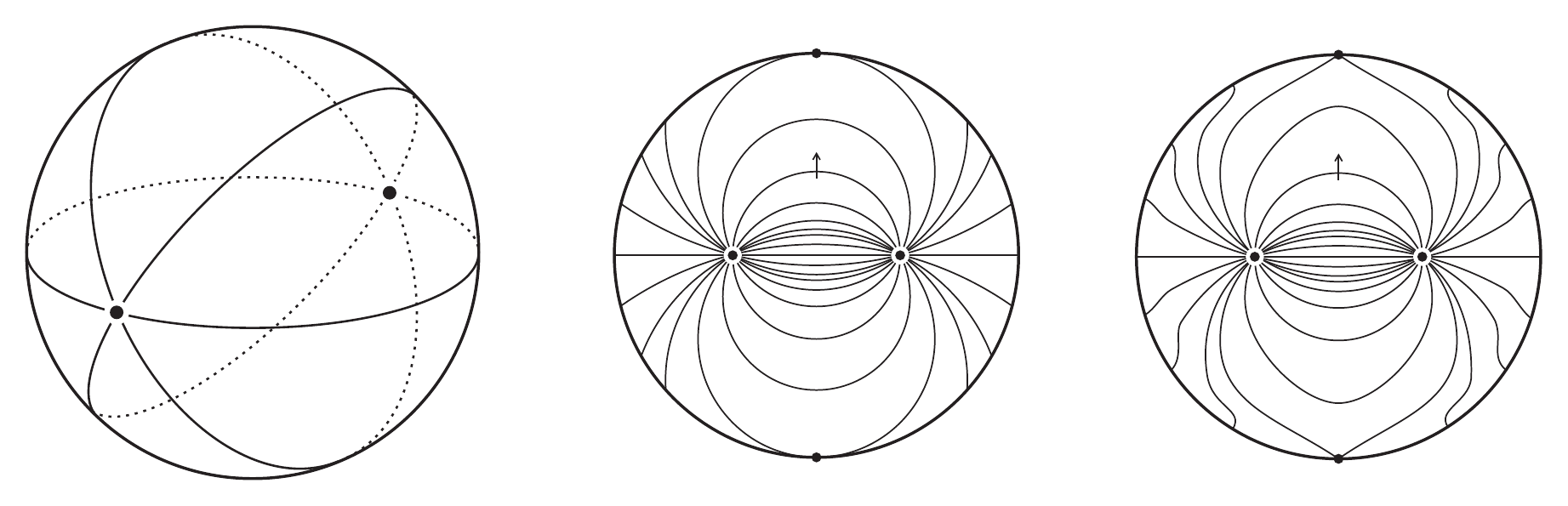}
\caption{}
\end{center}
\end{figure}

Let $D_-$ be a small 2-disc around any regular 
point $a$ of $f$. Denote by $D_+$
the complement $S^2\sm \Int D_-$,
so that $S^2=D_+\cup D_-$ and the discs 
$D_\pm$ intersect by their common boundary $\Sigma$.
Removing $D_-$
we obtain a map
$f:D_+\sm \{a,b\}\to S^1$;
the structure of its level lines is
shown  on the figure 1 (middle).


The restriction
$f~|~\Sigma$
has two non-degenerate critical points:
$N$ and $S$. 
The vector $v$ in the figure depicts the gradient of the
map $f$. Applying the rotation 
construction to $K_0$  we obtain 
a trivial 2-component link $R(K^0)$ in $S^3$.
Let $F_0 : S^3\sm R(K_0)\to S^1$ be
the unique $S^1$-invariant  function such that 
$F_0~|~ D_+ = f$.
This function is continuous, but not smooth,
since its level surfaces have conical singularities 
in the points of $\Sigma$.
To repair this, we will modify the function
$f$ in a \nei~ of $\Sigma$ 
so that the level lines of the modified function
$g:D_+\sm \{a,b\}\to S^1$ are
as depicted on the figure 1 (right).


Each non-singular level line intersecting $\Sigma$ is orthogonal to 
$\Sigma$ at the intersection point.
Let $G_0 : S^3\sm R(K_0)\to S^1$ be
the unique $S^1$-invariant  function such that 
$G_0~|~ D_+ = g$.
Then $G_0$ is a $\smo$ function having two critical 
points $N$ and $S$.
Observe that the 
descending disc
of the critical point $S$ of the function
$G_0~|~\Sigma$ 
is in $\Sigma$,
therefore the descending discs of $G_0$ will have the same dimension 
1, and $\ind_{G_0} S=1$.
The same reasoning holds for the 
{\it ascending disc } of the critical points $N$,
therefore
$\ind_{G_0} N=2$.

\subsection{The general case}
\lb{su:gen_case}

Let $\Si$ be the unit sphere in $\rr^{n+2}$, that is,
$$
\Si=\{(x_0,\ldots , x_{n+1}) ~|~  x_0^2+ \ldots   +x_{n+1}^2=1\}.$$
Denote by $\Si'$ its intersection with the hyperplane
$x_{n+1}=0$.
Let $a=(0,\ldots , 0, 1)$;
for each point
$z\in \Si'$
denote by
$C(z)$ the great
circle through $a, -a, z$, and by $C'(z)$
the closed semicircle containing these three points.
The projection $p$
onto the $(n+1)$-th coordinate
gives the bijection of 
$C'(z)$
onto the closed interval 
$[-1,1]$; this bijection 
is a diffeomorphism when restricted to
$C'(z) \sm \{a, -a\}$.
Let $\b:[-1,1] \to [-1,1]$
be a diffeomorphism
\sut~ $\b(x)=x$ for $x$ in a \nei~ of $\pm 1$.
Then there is a unique diffeomorphism $\bar\b$ of $\Si$ onto itself
\sut~ for every $z$ the curve $C'(z)$ is
$\bar\beta$-invariant and
$p(\bar\b(v))=\b(p(v))$ for every $v$.
The diffeomorphism $\bar\b$ will be 
called {\it the sliding, associated to $\b$}.
Observe that every sliding is isotopic to the identity map.

Let $D_\r\sbs \Si$ be the geodesic disc 
of radius $\r$ 
centered in $-a$. 
Let 
$$
D_-=D_{\pi/2} = \{(x_0, \ldots , x_{n+1}) ~|~ x_{n+1} \leq 0\},$$
$$
D_+= \{(x_0, \ldots , x_{n+1}) ~|~ x_{n+1} \geq 0\}.
$$
Put $\Si_\r=\pr D_\r$.
Let $N(\Si_\r, \e)$ denote the geodesic tubular \nei~ of $\Si_\r$.
For a given $\r$ and $\e>0$ sufficiently small
there is a sliding $\s$ sending $D_\r$ to $D_-$ 
and sending each normal geodesic segment of 
length $2\e$ in $N(\Si_\r, \e)$ isometrically to the corresponding 
 normal geodesic segment 
 in $N(\Si, \e)$. We have therefore a commutative diagram
 
 $$\xymatrix{
 N(\Si_\r, \e) \ar[r]^\s  & N(\Si, \e) \\
  \Si_\r\times ]-\e, \e[ \ar[u]^\Phi \ar[r]^{\bar\s}  & \Si\times ]-\e, \e[ \ar[u]^\Psi
 }$$
 where the vertical arrows are \df s and 
$\bar\s(x,\tau)=(\s(x),\tau)$.

Let $K$ be an $(n-1)$-knot in $S^{n+1}$ and
$\phi : S^{n+1}\sm K \to S^1$
a Morse map. We can assume that 
\been\item
$K\sbs \Int D_+$.
\item
$-a\notin Crit f$, 
\item
the \sma~ $\phi^{-1}(\phi(-a))$ is tangent to the hyperplane  
defined by the equation $x_n=0$.
\enen

The restriction $\phi~|~\pr D_\r$ can be considered as a real-valued Morse map.
Choosing $\r$ sufficiently small
we can assume that 
$\phi~|~\Si_\r$ is a Morse map having one  maximum and one minimum.
Denote the function
$\phi\circ \Phi$ by 
$h: 
\Si_\r\times ]-\e, \e[
\to
\rr.$
For $\r$ sufficiently  small, this function has the following property:

\begin{equation}\lb{D}
{\rm If \ \ }
\frac {\pr {h}}{\pr t}(x,t) =0,
\ \ 
{\rm then \ \ 
}
\frac {\pr {h}}{\pr x}(x,t) \not= 0,
{\rm where \ \ 
}
x\in \Si_\r, \ t\in ]-\e, \e[.
\end{equation}

Consider the restriction of $\phi$ to the subset
$S^{n+1}\sm (K\cup D_\r)$.
Composing $\phi$ with $\s^{-1}$ 
 we obtain a function
$$
\phi_0: D_+\sm \s(K) \to S^1.$$
This is a Morse map which extends to  
a geodesic tubular \nei~ 
of $\Si=\pr D^+$, and can be considered as a 
real-valued Morse function
in this \nei.
The restriction $\phi_0 ~|~ \Si$ 
has  two critical points of indices $n$ and $0$.
Denote these critical points by $N$ and $S$,
so that $\ind_{\phi_0} N=n, \ \ \ind_{\phi_0} S=0.$
The function $h_0=\phi_0\circ \Psi$ 
has the  following property:
\begin{equation}\lb{DD}
{\rm If \ \ }
\frac {\pr {h_0}}{\pr t}(x,t) =0,
\ \ 
{\rm then \ \ 
}
\frac {\pr {h_0}}{\pr x}(x,t) \not= 0,
{\rm where \ \ 
}
x\in \Si, \ t\in ]-\e, \e[.
\end{equation}
Now we will modify the function $\phi_0$ nearby $\Si$.
Let $\l:[-\e,\e]\to\rr$
be a $\smo$ function \sut~
$\l(t)=|t|$ for $t$ in a \nei~ of $\{-\e, \e\}$ and 
$\l(t)=t^2$ for $|t|\leq \e/2$.
Define a function
$h_1$ by the following formula:
$$ h_1(x,t)=h_0(x, \l(t)),
$$
and define a function 

$$\phi_1 :
D_+\sm \s(K) \to S^1
$$
as follows:

\been\item
if $v\notin N(\Si, \e)$, put $\phi_1(v)=\phi_0(v).$
\item
if $v\in N(\Si, \e), \ \ v=\Psi(x,t) $
with 
$x\in \Si, \ t\in ]-\e, \e[$, \ 
put 
$\phi_1(v) = h_1(x,t).$
\enen

\bepr\lb{p:crit_phi1}
The function $\phi_1$ has two critical  points in $N(\Si, \e)$, namely 
$N$ and $S$. Their indices are equal, respectively, to $n$ and $1$. 
\enpr
\Prf
The partial derivatives of $h_1$
are equal to 
$\frac {\pr {h_0}}{\pr x}(x,\l(t)) (x, t)$
and
$\frac {\pr {h_0}}{\pr t}(x,\l(t))\cdot \l'(t) $.
For $t=0$ the second derivative equals $0$,
and 
$\frac {\pr {h_1}}{\pr x}(x,\l(t)) (x, 0)$
vanishes in $N$ and $S$.
If $t\not=0$,
then 
$\l'(t)\not=0$,
and for $(x,t)$ to be a critical point
of $\phi_1$ it is necessary that 
$\frac {\pr {h_0}}{\pr t}(x,\l(t)) (x, t)$
vanish, which implies that
$\frac {\pr {h_0}}{\pr t}(x,t) (x, t)\not=0$
(see the property \rrf{DD}). $\qs$

Now we are ready to construct a Morse function
on the complement to $R(K)$.
Observe that the knot $K$ is 
equivalent to the knot $\s(K)$.
By  a certain abuse of notation we will replace $\s(K)$ by $K$, so in particular,
$K\sbs \Int D_+$.
Add one more coordinate $x_{n+2}$ and consider the sphere 
$$
\Si=\{(x_0,\ldots , x_{n+2}) ~|~  x_0^2+ \ldots   +x_{n+2}^2=1\}.$$
We have $D_+\sbs S^{n+2}$.
The knot $R(K)$ is defined by the following formula:

$$
R(K)
=
\{(x_0,\ldots , x_{n+2}) ~|~  
\Big(x_0,\ldots , x_{n}, \sqrt{x_{n+1}^2+x_{n+2}^2}\ \Big )\sbs K\}.$$
The circle $S^1$ acts on $S^{n+2}$ by rotations in the two last coordinates.
Define the Morse function
$\phi_2$ on the complement to
$R(K)$
by the two following properties:
\been\item
$\phi_2~|~D_+\sm K=\phi_1$.
\item $\phi_2$ is $S^1$-invariant.
\enen
The second property
implies that 
$$
\phi_2(x_0,\ldots , x_{n+2})
=
\phi_1\Big(x_0,\ldots , x_{n}, \sqrt{x_{n+1}^2+x_{n+2}^2}\ \  \Big ).
$$
Observe that the property 2) of the function $\phi_1$ guarantees that 
$\phi_2$ is $\smo$ on the subset $S^{n+2}\sm R(K)$.

\bepr\lb{p:crit_phi2}
\been\item
$Crit(\phi_2)= S^1\cdot Crit(\phi_1) \cup \{N,S\}$,
\item
The critical points $N$ and $S$ are 
non-degenerate, and
$$\ind_{\phi_2} N 
=
\ind_{\phi_1} N = n,
\ \ 
\ind_{\phi_2} S 
=
\ind_{\phi_1} S + 1 =2. $$
\enen 
\enpr 
\Prf
The point 1) is easy to deduce from the definition of $\phi_2$.
As for the indices of the critical points
observe that the descending disc of the critical point $N$ 
in $N(\Si, \e)$
belongs to the sphere $\Si$ 
which is fixed by the action of $S^1$.
Thus the index of $N$ does not change when we replace $\phi_1$ by $\phi_2$.
A similar argument applies to the {\it ascending disc } of $S$,
and this implies the rest of the proposition. $\qs$

Each critical point of $\phi_1$ gives rise to a circle of 
critical poins of $\phi_2$. Using the same method, as in 
the previous work of the authors,
we 
wwe perturb the function $\phi_2$ in a neighbourhood of 
each of these critical circles, and obtain finally a \rmf~ $R(\phi)$
on the complement to $R(K)$
\sut~
$$
\# Crit(R(\phi))= 2 \# Crit(\phi_2) + 2.
$$
This completes the proof of Theorem 
\ref{t:mn_rotation}.
$\qs$ 

\subsection{4-thread spinning}
\lb{su:4thread}

In this subsection we give a brief description of 
one more construction of surface-links.
Let $L\sbs S^3$ be a classical link 
and $\phi:S^3\sm L\to S^1$ a Morse map.
Let $p,q\in L$ and let $\g:[0,1]\to S^3$ be a $\smo$ curve joining 
$p$ and $q$ and belonging entirely to one of the regular level surfaces 
$\phi^{-1}(\l)$ of the map $\phi$.
We assume morover that $\Im\g \cap L=\{p,q\}$, and 
that $\g'(0)$ and $\g'(1)$ are not tangent to
$L$. 
Let $D$ be a small \nei~ of $\Im\g$ diffeomorphic to
a 3-disc. Denote by $\Sigma$ its boundary.
We can assume that $L\cap \Si$
consists of four points and that $L$ 
is orthogonal to $\Sigma$ 
at each of these points.
Denote by $S_0^2$ the 2-sphere with 4 points removed.
Recall that there is a standard Morse function
$\phi_0$ on $S^2_0$ having 2 critical points of indices 1.
We can assume  that 
the restriction of $\phi$ to $\Si\sm L$
is diffeomorphic to $\phi_0$.

Remove the interior of $D$ from $S^3$
and rotate the remaining manifold $S^3\sm \Int D$ 
around $\Si$. We obtain the sphere $S^4$; 
the subset which is spun by $L\sm \Int D$ 
during the rotation is an embedded 2-surface in $S^4$.

We call this construction
{\it 4-thread spinning}
to distinguish it from the usual spinning,
and denote the resulting surface-link by $S'(L)$.
If $p$ and $q$ are on different 
connected components 
of $L$,
then the number of
connected components 
of 
 $S'(L)$
 is the same as for $L$.
 If $p$ and $q$ are in different
  connected components 
of $L$, then the number of
connected components 
of $S'(L)$ equals that of $L$ increased by $1$.
Applying the same method as in the Subsection
\ref{su:gen_case}
we can construct a Morse function 
$\wi \phi$ on $S^4\sm S'(L)\to S^1$
such that 
$m(\wi\phi)= m(\phi)+2.$
\beco\lb{co:4thread}
$$
\MM\NN(S'(L))\leq 2\MM\NN(L)+2.\hspace{3cm}\qs$$
\enco

\section{Surface-links}
\lb{s:surf-links}


In this section we develop 
circle-valued Morse theory for 
surface-links.

\subsection{Motion pictures and saddle numbers}
\lb{su:saddle}

Let $F$ be a surface-link,
that is, a closed oriented 2-dimensional 
$\smo$ submanifold of $S^4$. We can assume $F\sbs \rr^4$.

Choose a  projection $p$
of $\rr^4$ onto a line. 
Assume that the critical points of the 
function $p|F$ are non-degenerate.
Denote by $sdl(F)$ the minimal number of saddle points of $p|F$
over all the projections $p$.

\bede\lb{d:saddle}
A saddle number $sd(F)$ is the minimum of numbers $sdl(F')$
where $F'$ ranges over all surface-links $F'$ ambiently isotopic to 
$F$.
\end{defi}

The invariant $sd(F)$ is closely related to the {\it $ch$-index }
of $F$, introduced and studied by K. Yoshikawa in \cite{Yoshikawa}.
In particular, we have $sd(F)\leq ch(F)$. 
In order to relate the  number $sd(F)$ to $\mnk$ we 
will reformulate the definition of the saddle number.

Let $F\subset S^4$ be a 
surface-link.
The equatorial $3$-sphere $\Sigma^3$ of the standard Euclidean sphere $S^4$ 
divides $S^4$ into two parts: 
$$
S^4=D^4_+\cup D^4_-,\ \  {\rm with }\ \  D^4_+\cap D^4_-=\Sigma^3.
$$
We assume that $F$ is included in ${\rm Int}(D^4_-)$ and 
$F$ does not include the centre of $D^4_-$. 
Perturbing the embedding $F\subset D^4_-$ if necessary, 
we can assume that the restriction $\rho=r|_F$ 
of the radius function $r:D^4_-\rightarrow [0,1]$ 
is a Morse function. 
The family $\{(r^{-1}(t), \rho^{-1}(t))\}_{t\in [0,1]}$ of possibly
singular links can be drawn 
as a {\it motion picture} (see \cite{Kamada2002}, Chapter 8). 
Each singularity of a link in the family corresponds to a critical point of $\rho$. 
A critical point of $\rho$ of index $0$ ($1$, $2$, respectively) is called 
{\it minimal point} ({\it saddle point}, {\it maximal point}, respectively) of $\rho$, 
which is represented by a {\it minimal band} ({\it saddle band}, 
{\it maximal band}, respectively) 
in (a modification of) the motion picture. 

It is clear that 
the minimal number of the saddle points for all such Morse functions $\rho$ 
is equal to $sd(F)$.

\beth\lb{t:mn-s} \ \ 
$\mathcal{MN}(F)\leq 2\, sd(F)+\chi(F)-2$. 
\enth
{\it Proof}. 
Since $\rho$ is a Morse function, 
the manifold $D^4_-\setminus {\rm Int}\, N(F)$ admits a handle decomposition 
with 
{\bc
one $0$-handle and
}
$m_i(\rho)$ $(i+1)$-handles for $i\in\{0,1,2\}$ (see 
{\bc
\cite{G}, and  also}
\cite{GS1999}, Proposition 6.2.1). 

The exterior $E(F)=S^4\setminus {\rm Int}\, N(F)$ of $F$ is obtained by attaching 
a $4$-handle $D^4_+$ to $D^4_-\setminus {\rm Int}\, N(F)$. 
Since $D^4_-\setminus {\rm Int}\, N(F)$ is connected, 
there is a $3$-handle in $D^4_-\setminus {\rm Int}\, N(F)$ 
which connects $\partial N(F)$ with $\partial D^4_-$. 
Thus the $3$-handle cancels the $4$-handle $D^4_+$ (see \cite{Milnor1965}, Section 5). 
Turning the handlebody upside down, 
we obtain a dual decomposition of $E(F)$ and 
a corresponding Morse function $f:E(F)\rightarrow \mathbb{R}$ 
which is constant on $\partial E(F)$ and the following Morse numbers: 
$m_1(f)=m_2(\rho)-1$, $m_2(f)=m_1(\rho)$, $m_3(f)=m_0(\rho)$, $m_4(f)=1$. 

Using the argument from work of the second author \cite{P-t}, p. 629,
we can deform the real-valued Morse function $f$
to a circle-valued regular function $\phi:E(F)\to S^1$,
such that $m_k(f)=m_k(\phi)$ for every $k$.
Consider the function $-\phi$, which has one critical point 
of index $0$. Applying the cancellation of this local minimum,
we obtain a Morse function $\psi:E(F)\to S^1$
belonging to the class $-\xi$, and such that 
$m_0(\psi)=0,\ m_1(\psi)=m_3(f)-1, \ m_2(\psi)=m_2(f), \ m_3(\psi)=m_1(f), \ m_4(\psi)=0$.
Put $g=-\psi$. Then we have
$$
m_0(g)=m_4(g)=0, \ \ 
m_1(g) = m_2(\r) -1, \ $$
$$
m_2(g) = m_1(\r), \ \ 
m_3(g) = m_0(\r) -1.
$$
Observe that $m_0(\r)-m_1(\r)+m_2(\r)=\chi(S^2)=2$, 
therefore the total number of critical points of 
$g$ equals $2m_1(\r)$. Choosing the function $\r$ with 
$m_1(\r)=sd(F)$ we accomplish the proof. $\qs$


\beco\lb{c:2knots}
Let $K\sbs S^4$ be a 2-knot. Then 
$\MM\NN(C_K)\leq 2sd(K).\qs$
\enco

\bepr\lb{p:trivial}
Let $F\sbs S^4$ be the trivial $k$-component surface-link.
Then 

$\MM\NN(F)=4k-2-\chi(F).$
\enpr\Prf
It is not diffcult to show that 
$\wh b_1(C(F))\geq k-1, \ \wh b_3(C(F))\geq k-1$.
Therefore for every regular Morse map
$f:C(F)\to S^1$ we have $m_1(f)+m_3(f)\geq 2(k-1)$. Assuming $m_0(f)=m_4(f)=0$
we have $m_1(f)-m_2(f)+m_3(f) =2-\chi(F)$, and $\MM\NN(C(F))\geq 4k-2-\chi(F)$;
this lower bound coincides with the upper bound derived from
Theorem \ref{t:mn-s}. $\qs$

\subsection{Spun knots}
\lb{su:spun}

Let $K$ be a classical knot in $S^3$
denote by $S(K)$ the corresponding spun knot. 

\bepr\lb{p:mn-spun}
If $K$ is a non-fibered knot of tunnel number 1, then
$\MM\NN(S^4\sm S(K))=4$.
\enpr
\Prf
Recall that 
$\MM\NN(S^4\sm S(K)) \leq 2\MM\NN(K)$
(Corollary \ref{c:class-spin}).
In the paper \cite{P-t} of the second author it is shown that 
$\MM\NN(K)\leq 2t(K)$, hence $\MM\NN(S(K))\leq 4$
by Corollary \ref{c:class-spin}.
Put $G=\pi_1(S^3\sm K)$, then $\pi_1(S^4\sm S(K))$;
let $H=[G,G]$.
Let $f: S^4\sm S(K) \to S^1$ be a regular Morse map without minima and maxima.
If $m_1(f)=0$, then a standard Morse-theoretic argument 
applied to the infinite cyclic cover of $\ssmk$ 
implies that $H$ 
is finitely generated, which is impossible, since $K$ is not fibred.
Therefore $m_1(f)\geq 1$, and similarly, $m_3(f)\geq 1$,
hence $m_2(f)\geq 2$ and the proposition is proved. $\qs$

\subsection{Surface-links of Yoshikawa's table}
\lb{su:yosh_links}

Yoshikawa  \cite{Yoshikawa}
suggested a method for enumerating 
surface-links.
To each surface-link $F$
he associated a natural number
$ch(F)$. His methods allowed him to make a list 
of all (weakly prime) 
surface-links $F$ with $ch(F)\leq 10$.
It is clear from the definition
of the invariant $ch(F)$ that we have 
$sd(F)\leq ch(F)$.
In the rest of this section we assume that the reader
is familiar with Yoshikawa's work,
and with his terminology.
There are 6 two-knots
in Yoshikawa's table,
namely

$$0_1, \ 8_1, 9_1, \ 10_1, \ 10_2, \ 10_3.$$
The trivial 2-knot $0_1$ is obviously fibred.
The knots $8_1$ and $10_1$ are spun knots
of the trefoil knot and respectively of the figure 8
knot, thus both $8_1$ and $10_1$ are fibred by \cite{AS}.

The case of $9_1$ is more complicated.
The saddle number of this 2-knot is 2.
Therefore $\MM\NN(9_1)\leq 4$.
Using the presentation of the fundamental group
 of the complement to $9_1$ (see \cite{Yoshikawa})
and Poincar\'e duality properties
it is easy to compute the Novikov numbers 
of $9_1$. Namely we have
$\wh q_1=1, \wh q_2=\wh q_3=0.$
Therefore 
$$
2\leq \MM\NN(9_1)\leq 4.
$$
 The 2-knot $10_2$ 
is the 2-twist-spun knot of the trefoil knot,
hence fibered by Zeeman's theorem \cite{Z}.
Similarly, $10_3$ is fibered, being the 3-twist spun of 
the trefoil knot. 

The surface-link $6_1^{0,1}$
is the result of spinning of the Hopf link which is fibred
(see the left of Figure 2)
therefore $\MM\NN (6_1^{0,1}) = 0.$

The surface-link $8_1^{1,1}$ is the spun torus of 
the Hopf link. Applying Theorem \ref{t:mn_rotation}
we get the upper bound $\MM\NN (8_1^{1,1}) \leq 2.$ 
Computing the Euler charcateristic implis 
the inverse inequality, so $\MM\NN (8_1^{1,1}) = 2.$

The same argument applies to the surface-link
$10_1^1$, which is the spun torus of the trefoil knot, see the figure 2 (middle), 
so that $\MM\NN (10_1^1) = 2.$

The surface-link $10_1^{0,1}$
is the result of spinning of the link $4^2_1$ which is fibred,
therefore $\MM\NN (10_1^{0,1}) = 0.$

The case of the surface-link $F=10_1^{0,0,1}$ is more complicated.
This surface-link is the result of 4-threaded spinning 
of the connected sum $L$ of two copies of the Hopf link,
see Figure 2 (right)
and applying Corollary \ref{co:4thread}
we deduce 
$\MM\NN (F)\leq 2$.
The computation of Euler characteristic 
gives the lower bound 2 for the Morse-Novikov
number, thus
$\MM\NN (10_1^{0,0,1})=2$.

\begin{figure}
\begin{center}
\includegraphics{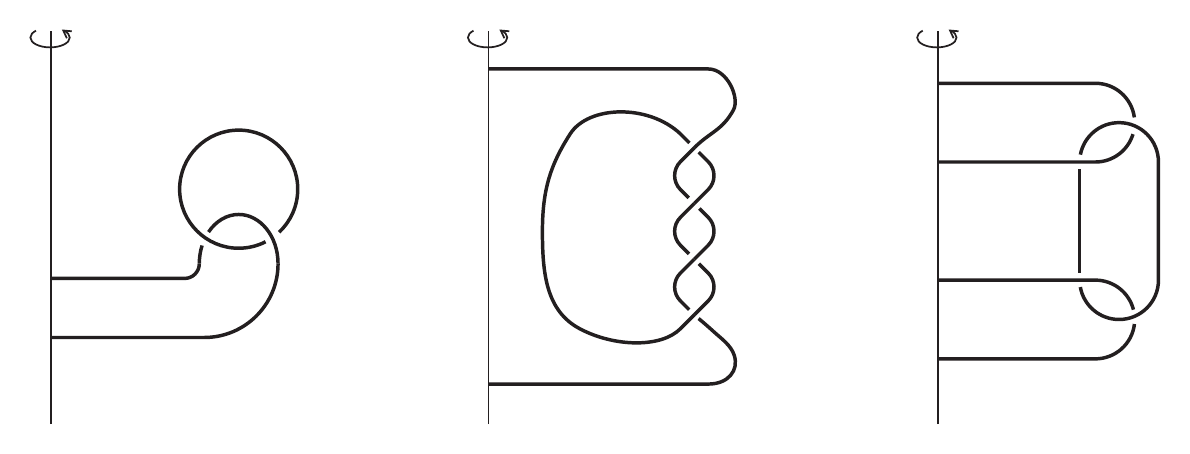}
\caption{}
\end{center}
\end{figure}

\section{Acknowledgements}
\lb{s:ac}

This work was accomplished when the second author 
was visiting the Tokyo Institute of Technology in 2016
with the support of the JSPS fellowship.
The first author was partially supported by JSPS KAKENHI
Grant Numbers 25400082, 16K05142.
The second author thanks the Tokyo Institute of Technology
for support and  warm hospitality.


\end{document}